\RequirePackage{fix-cm}
\documentclass[smallcondensed]{svjour3}     
\smartqed  
\usepackage{lipsum}
\usepackage{amsfonts}
\usepackage{graphicx}
\usepackage{epstopdf}
\usepackage{graphicx}
\usepackage{color,latexsym}
\usepackage{verbatim}
\usepackage{lscape}
\usepackage{amsmath}
\usepackage{amssymb}
\usepackage{esint}
\usepackage{multirow}
\usepackage{hhline}
\usepackage{tabularx}
\usepackage{mathtools}
\usepackage{float}
\usepackage{setspace}

\newcommand{\defeq}{\vcentcolon=}

\newcommand{\sign}{\operatorname{sign}}
\ifpdf
  \DeclareGraphicsExtensions{.eps,.pdf,.png,.jpg}
\else
  \DeclareGraphicsExtensions{.eps}
\fi

\begin{document}

\title{Approximate Integrals Over Bounded Volumes with Smooth Boundaries\thanks{}
}


\author{Jonah A. Reeger}


\institute{J. A. Reeger \at
Department of Mathematics and Statistics\\ Air Force Institute of Technology\\2950 Hobson Way\\Wright-Patterson Air Force Base, Ohio\\
              United States \\
              \email{jonah.reeger@afit.edu} }

\date{Received: date / Accepted: date}

\maketitle

\begin{abstract}
A Radial Basis Function Generated Finite-Differences (RBF-FD) inspired technique for evaluating definite integrals over bounded volumes that have smooth boundaries in three dimensions is described.  A key aspect of this approach is that it allows the user to approximate the value of the integral without explicit knowledge of an expression for the boundary surface.  Instead, a tesselation of the node set is utilized to inform the algorithm of the domain geometry.  Further, the method applies to node sets featuring spatially varying density, facilitating its use in Applied Mathematics, Mathematical Physics and myriad other application areas, where the locations of the nodes might be fixed by experiment or previous simulation. By using the RBF-FD-like approach, the proposed algorithm computes quadrature weights for $N$ arbitrarily scattered nodes in only $O(N\mbox{ log}N)$ operations with tunable orders of accuracy.
\keywords{Radial Basis Function \and RBF \and quadrature \and volume}
 \subclass{68Q25 \and  65R99}
\end{abstract}

\section{Introduction}
This article is concerned with the development of a method for the approximate evaluation of definite integrals over bounded volumes that have smooth boundaries in $\mathbb{R}^{3}$.  That is, consider volume integration over the domain, $\Omega$, where the boundary, $\partial\Omega$, is described by, for instance, all points satisfying $h(\mathbf{x})\leq 0$. Here $h:\mathbb{R}^{3}\mapsto\mathbb{R}$ is assumed to be a smooth function.  A key aspect of the algorithm presented here, however, is that the mapping $h$ need not be known.  All that is required is a tesselation of the domain with tetrahedra.

This common problem of approximating the values of definite integrals, where some of the earliest work traces back many centuries to early attempts to measure the area of the circle \cite{WFMG2018}, is known as quadrature when considering integration over an interval, or quadrature/cubature when considering integration occurs over domains in two or more dimensions. Now there are numerous texts devoted to summarizing sophisticated and accurate techniques for estimating the values of integrals over intervals, areas, and volumes (see, for example \cite{WFMG2018,AHS1971,ARKCWU1998,PKKMRS2005}).

The problem of evaluating definite integrals over volumes with smooth bounding surfaces is gaining popularity.  In particular, attention has recently been paid to improvements to finite element methods in the presence of implicitly defined surfaces, but often with illustrations only in $\mathbb{R}^{2}$ \cite{MAODS2016} or for elements which impose a particular structure on the node sets, e.g. hyperrectangles \cite{RIS2015}.  These methods can often rely on choosing appropriate coordinate systems for applying Gaussian quadrature over each dimension (see, e.g. \cite{TCWLHLLZWZ2020}). Further, some of these previous methods required explicit knowledge of $h$ and its derivatives, an unnecessary complication remedied in section \ref{sec:Unknown_Param}.  Other methods for integrating over volumes rely on the Gauss divergence theorem \cite{BMFKMO2013} and node sets designed for certain characteristics \cite{VKRMKAN2018}.

Since numerical quadrature is often a follow-up to some other task (such as collecting data or numerically solving Partial Differential Equations (PDEs)), it can be impractical to require node locations that are specific to the quadrature method. It is true that an extra interpolation step could take place to approximate the integrand at node locations specific to the quadrature rule, but this is an unnecessary step that imposes constraints on the compatibility of the accuracy of the interpolation and quadrature. Further, such a step may be unworkable if the interpolation step itself places constraints on the node set. Requirements on node locations are common in, for instance, Newton-Cotes or Gaussian Quadrature rules for evaluating definite integrals of intervals in 1-dimension.   In both cases a polynomial interpolant of the integrand is constructed and then integrated, with interpolation conditions enforced at node locations that either have equal spacing between adjacent nodes or are tied to roots of orthogonal polynomials.

To overcome these requirements on the structure of the node set, the present algorithm is designed to find the quadrature weights given a node set defined by the application or user.  It is possible to construct a polynomial interpolant of the integrand in this context, followed by integration of the interpolant, in more than one dimension.  Unfortunately, however, the polynomial basis set used for interpolation does not depend on the locations of the nodes, and such basis sets suffer from a question of the existence and uniqueness of an interpolant on unstructured node sets \cite{PCC1959,JCM1956}, which has prompted the use of Radial Basis Functions (RBFs) in the approximation of the integrand.

This paper describes an extension of the work in \cite{JAR2020}, which considered only the volume of the ball in $\mathbb{R}^{3}$, to volumes bounded by arbitrary smooth surfaces. In earlier works by the author and collaborators the concept of RBF-FD (radial basis function-generated finite differences), which was designed for and has been successfully applied to solve PDEs (see \cite{BFNF2015a,BFNF2015b} for surveys of these applications), was generalized to construct quadrature rules for an interval in one-dimension, bounded domains in two-dimensions and, more specifically, integrals over bounded two-dimensional (piecewise-)smooth surfaces embedded in three-dimensions (see \cite{JARBF2016,JARBFMLW2016,JARBF2017} and the references therein).

The following Section \ref{sec:Algorithm_Description} describes the present quadrature method.  Section \ref{sec:Test_Examples} describes some test examples, with illustrations of convergence rates and computational costs. Finally, section \ref{sec:Conclusions} outlines some conclusions. An implementation of the methods discussed in this paper has been made available at \cite{VolumeQuadCodeMatlab}.

\section{Description of the key steps in the algorithm} \label{sec:Algorithm_Description}

Consider evaluating
\begin{align}
\iiint\limits_{\Omega}f(\mathbf{x})dV, \label{eq:volume_integral}
\end{align}
where $\Omega=\left\{\mathbf{x}\in\mathbb{R}^{3}:h(\mathbf{x})\leq0\right\}$.  Here, the boundary of $\Omega$ can be defined as $\partial \Omega = \left\{\mathbf{x}\in\mathbb{R}^{3}:h(\mathbf{x})=0\right\}$.

Similar to the work presented in \cite{JARBF2016,JARBFMLW2016,JARBF2017} for surface integrals and \cite{JAR2020} for the volume of the ball in $\mathbb{R}^{3}$, at a high level, the proposed algorithm can be described in four steps:
\begin{enumerate}
    \item Decompose the domain of integration into $K\in \mathbb{Z}^{+}$ subdomains.
    \item On the $k^{\mbox{th}}$ subdomain define an interpolant of the integrand.
    \item Integrate the basis functions defining the interpolant and solve a system of linear equations to determine weights for integrating an arbitrary function $f$ over the $k^{\mbox{th}}$ subdomain.
    \item Combine the weights for the integrals over the $K$ subdomains to obtain a weight set for approximating the volume integral of $f$ over $\Omega$.
\end{enumerate}
Each of these steps is described in greater detail in what follows.  Throughout the remainder of this paper, the subscript $k$ will be used to indicate that the steps are carried out for each tetrahedron separately.  When further subscripting is necessary--for instance, to indicate entries of a vector, matrix, or a set--the necessary indexing will follow after a comma.  

\subsection{Step 1: Decompose the domain of integration} \label{sec:Step1}

Suppose that $\mathcal{S}_{N}=\left\{\mathbf{x}_{i}\right\}_{i=1}^{N}$ is a set of $N$ unique points in $\Omega$, with a subset exactly on the boundary surface. On the set $\mathcal{S}_{N}$ construct a tessellation $T=\left\{t_{k}\right\}_{k=1}^{K}$ (via Delaunay tessellation or some other algorithm) of $K$ tetrahedra.  These tetrahedra encompass the bulk of the volume of $\Omega$.  However, the tetrahedra near the boundary, with $h=0$ on at least three vertices, likely over- or under-approximate the volume locally.

Let $\mathcal{K}_{S}\subset\{1,2,\ldots,K\}$ be the set of indices such that if $k\in \mathcal{K}_{S}$ then the tetrahedron $t_{k}$ has at least one face that is not shared with any of the other tetrahedra, call it $\tau_{k,*}$.  This face has all three vertices on the surface of $\Omega$, i.e. $h=0$ on these vertices, and unless the surface is planar there is a sliver of volume, $s_{k}$, between $\tau_{k,*}$ and the smooth, curved bounding surface that must be accounted for in decomposing the volume.  Conversely, let $\mathcal{K}_{I}= \left\{1,2,\ldots,K\right\}\backslash \mathcal{K}_{S}$ be the indices of tetrahedra that do not have a face with three vertices on the surface.  With these definitions \eqref{eq:volume_integral} can be decomposed as
\begin{align}
\iiint\limits_{\Omega}f(\mathbf{x})dV=\sum\limits_{k\in \mathcal{K}_{I}}\iiint\limits_{t_{k}}f(\mathbf{x})dV+\sum\limits_{k\in \mathcal{K}_{S}}\left(\iiint\limits_{t_{k}}f(\mathbf{x})dV+\nu_{k}\iiint\limits_{s_{k}}f(\mathbf{x})dV\right), \label{eq:volume_integral_decomposed}
\end{align}
with $\nu_{k}^2=1$ and $\nu_{k}\in\mathbb{R}$.  The region $s_{k}$ can exist entirely inside, entirely outside, or might be partially inside and outside of $t_{k}$.  Therefore, to make \eqref{eq:volume_integral_decomposed} valid, it is important to choose a coordinate system and $\nu_{k}$ so that only portions of $s_{k}$ outside of $t_{k}$ are added and the remainder is subtracted.  A useful way of defining the coordinate system is to locate the origin in the plane containing $\tau_{k,*}$ and to align two of the coordinate axes so that they are parallel to $\tau_{k,*}$.  Two such coordinate systems will be defined in section \ref{sec:Sliver_Integrals}.

\subsection{Step 2: Define an interpolant of the integrand \label{sec:ConstructInterpolant}}

For each tetrahedron in $T$ define the sets $\mathcal{N}_{k}=\left\{\mathbf{x}_{k,j}\right\}_{j=1}^{n}$ to be the $n$ points in $\mathcal{S}_{N}$ nearest to the midpoint, $m_{k}$, of $t_{k}$ (the average of the vertices of $t_{k}$).  Then for each $k$ the integrals
\begin{align}
\iiint\limits_{t_{k}}f(\mathbf{x})dV\nonumber
\end{align}
and
\begin{align}
\iiint\limits_{s_{k}}f(\mathbf{x})dV\nonumber
\end{align}
are evaluated by first approximating $f(\mathbf{x})$ by an RBF interpolant, with interpolation points from the set $\mathcal{N}_{k}$, and then integrating the interpolant. In this work, the interpolant can be written as a linear combination of (conditionally-) positive definite RBFs,
\begin{align}
\phi\left(\left\lVert \mathbf{x}-\mathbf{x}_{k,j}\right\rVert_{2}\right), j=1,2,\ldots,n\nonumber
\end{align}
and (supplemental) multivariate polynomial terms.  If $k\in\mathcal{K}_{S}$, then the same interpolant and set of nodes is used for approximating the integrand over both $t_{k}$ and $s_{k}$.  Define $\{\pi_{k,l}(\mathbf{x})\}_{l=1}^{M}$, with $M=\frac{(m+1)(m+2)(m+3)}{6}$, to be a set of all of the trivariate polynomial terms up to degree $m$. The interpolant is constructed as
\begin{align}
s(\mathbf{x}):=\sum_{j=1}^{n}c_{k,j}^{\mbox{RBF}}\phi\left(\left\lVert \mathbf{x}-\mathbf{x}_{k,j}\right\rVert_{2}\right)+\sum_{l=1}^{M}c_{k,l}^{p}\pi_{k,l}(\mathbf{x})=\mathbf{c}_{k}^{T}\mathbf{\Phi}_{k}(\mathbf{x}),\label{eq:rbf_interpolant}
\end{align}
where the entries of the $(n+M)\times 1$ coefficient vector
\begin{align}
\mathbf{c}_{k}=\left[\begin{array}{cccccccc} c_{k,1}^{RBF} & c_{k,2}^{RBF} & \cdots & c_{k,n}^{RBF} & c_{k,1}^{p} & c_{k,2}^{p} & \cdots & c_{k,M}^{p}\end{array}\right]^{T}\nonumber
\end{align}
are chosen to satisfy the interpolation conditions $s(\mathbf{x}_{k,j})=f(\mathbf{x}_{k,j})$, $j=1,2,\ldots,n$, along with the typical constraints applied to RBF interpolants, $\sum_{j=1}^{n}c_{k,j}^{RBF}\pi_{l}(\mathbf{x}_{k,j})=0$, $l=1,2,\ldots,M$. The $(n+M)\times 1$ vector $\mathbf{\Phi}_{k}(\mathbf{x})$ consists of all of the basis functions evaluated at $\mathbf{x}$, i.e.
\begin{align}
    \lbrack\mathbf{\Phi}_{k}(\mathbf{x})\rbrack_{j}=\left\lbrace\begin{array}{cc} \phi\left(\left\lVert \mathbf{x}-\mathbf{x}_{k,j}\right\rVert_{2}\right) & j=1,2,\ldots,n \\ \pi_{k,j-n}(\mathbf{x}) & j=n+1,\ldots,n+M\end{array}
    \right..\nonumber
\end{align}

Satisfaction of the interpolation conditions and constraints amounts to solving the system of linear equations $A_{k}\mathbf{c}_{k}=\mathbf{f}_{k}$ with the $(n+M)\times(n+M)$ matrix
\begin{align}
A_{k}=\left[\begin{array}{cc}\Phi_{k}^{T} & P_{k}\\P_{k}^{T} & 0\end{array}\right].\nonumber
\end{align}
The $n\times n$ submatrix $\Phi_{k}$ is made up of the RBFs evaluated at each point in $\mathcal{N}_{k}$, that is
\begin{align}
\Phi_{k,ij}=\phi\left(\left\lVert \mathbf{x}_{k,i}-\mathbf{x}_{k,j}\right\rVert\right),\mbox{ for }i,j=1,2,\ldots,n.\nonumber
\end{align}
Likewise the $n\times M$ matrix $P_{k,il}$ consists of the polynomial basis evaluated at each point in $\mathcal{N}_{k}$ so that
\begin{align}
P_{k,il}=\pi_{k,l}(\mathbf{x}_{k,i}),\mbox{ for }i=1,2,\ldots,n\mbox{ and }l=1,2,\ldots,M.\nonumber
\end{align}
The right hand side consists of the function $f$ evaluated at the points in $\mathcal{N}_{k}$ and an $M\times 1$ vector of zeros, $\mathbf{0}_{M}$, that is,
\begin{align}
\mathbf{f}_{k}=\left[\begin{array}{ccccc}f(x_{k,1}) & f(x_{k,2}) & \cdots & f(x_{k,n}) & \mathbf{0}_{M}^{T} \end{array}\right]^T,\nonumber
\end{align}
In the context of constructing a quadrature rule that is applicable to more than one integrand at a time, i.e. to avoid computing the values of $c_{k,j}^{\mbox{RBF}}$ and $c_{k,l}^{p}$, it is convenient to write the interpolant as (noting that $\mathbf{c}_{k}=A_{k}^{-1}\mathbf{f}_{k}$)
\begin{align}
    s(\mathbf{x})= (A_{k}^{-1}\mathbf{\Phi}_{k}(\mathbf{x}))^{T}\mathbf{f}_{k}.\nonumber
\end{align}
The vector of functions $\mathbf{\Psi}_{k}(\mathbf{x})=A_{k}^{-1}\mathbf{\Phi}_{k}(\mathbf{x})$ is a cardinal basis for the space spanned by the functions in $\mathbf{\Phi}_{k}(\mathbf{x})$.

\subsection{Step 3: Integrate the interpolant of the integrand \label{sec:IntegrateInterpolant}}
By integrating the interpolant, the approximation of the integral of $f$ is reduced to
\begin{align}
\iiint\limits_{t_{k}}f(\mathbf{x})dV\approx\iiint\limits_{t_{k}}s(\mathbf{x})dV=\sum_{j=1}^{n}w_{k,j}f(\mathbf{x}_{k,j})\nonumber
\end{align}
for $k\in\mathcal{K}_{I}$ and
\begin{align}
\iiint\limits_{t_{k}}f(\mathbf{x})dV+\nu_{k}\iiint\limits_{s_{k}}f(\mathbf{x})dV\approx\iiint\limits_{t_{k}}s(\mathbf{x})dV+\nu_{k}\iiint\limits_{s_{k}}s(\mathbf{x})dV=\sum_{j=1}^{n}w_{k,j}f(\mathbf{x}_{k,j})\nonumber
\end{align}
for $k\in\mathcal{K}_{S}$.  Given the definition of $\mathbf{f}_{k}$, the quadrature weights, $w_{k,j}$, are the first $n$ entries of the solution to
\begin{align}
    A_{k}\mathbf{w}_{k}=\mathbf{I}_{k}.\nonumber
\end{align}

If $k\in\mathcal{K}_{I}$, the right hand side, $\mathbf{I}_{k}$, includes integrals of the basis functions over $t_{k}$ only.  That is,
\begin{align}
I_{k,j}=\left\{\begin{array}{cc} \iiint\limits_{t_{k}}\phi\left(\left\lVert \mathbf{x}-\mathbf{x}_{k,j}\right\rVert\right)dV & j=1,2,\ldots,n \\ \iiint\limits_{t_{k}}\pi_{j-n}(\mathbf{x})dV & j=n+1,n+2,\ldots,n+M\end{array}\right..\nonumber
\end{align}
The integrals of the trivariate polynomial terms can be evaluated exactly via, for instance, the Divergence Theorem or barycentric coordinates.  For the RBFs, the integrals can be evaluated by further decomposing $t_{k}$ into four tetrahedra that share the common vertex $\mathbf{x}_{k,j}$.  Summing the integrals of the RBFs over the four tetrahedra results in the integral over $t_{k}$.  This process allows the volume integral over a tetrahedron to be reduced to four integrals in a single dimension.  Section 2.3.1 of \cite{JAR2020} explains this process.

On the other hand, for $k\in\mathcal{K}_{S}$
\begin{align}
I_{k,j}=\left\{\begin{array}{cc} \iiint\limits_{t_{k}}\phi\left(\left\lVert \mathbf{x}-\mathbf{x}_{k,j}\right\rVert\right)dV+\nu_{k}\iiint\limits_{s_{k}}\phi\left(\left\lVert \mathbf{x}-\mathbf{x}_{k,j}\right\rVert\right)dV & j=1,2,\ldots,n \\ \iiint\limits_{t_{k}}\pi_{j-n}(\mathbf{x})dV+\nu_{k}\iiint\limits_{s_{k}}\pi_{j-n}(\mathbf{x})dV  & j=n+1,\ldots,n+M\end{array}\right..\nonumber
\end{align}
The integrals over $t_{k}$ are evaluated using the methods described in the previous paragraph and in section 2.3.1 of \cite{JAR2020} while the integrals over $s_{k}$ are approximated using a scheme discussed in section \ref{sec:Sliver_Integrals} of this paper.

\subsubsection{Integrals Over Slivers of Volume at the Surface} \label{sec:Sliver_Integrals}

\begin{figure}[h]
\begin{center}
\includegraphics[width=1.5in]{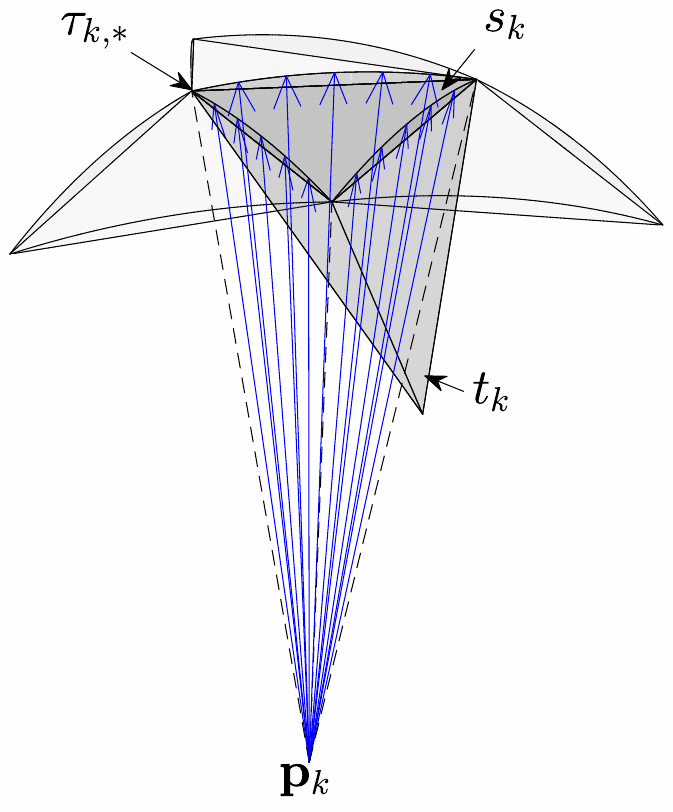}
\end{center}
\caption{An illustration of one of the tetrahedra in the set $T$ near the bounding surface.  Call the tetrahedron shaded in darker grey $t_{k}$ and let $\tau_{k,*}$ be the triangular face of $t_{k}$ with three vertices on the bounding surface.  The arrows originating at $\mathbf{p}_{k}$ indicate the projection of $\tau_{k,*}$'s edges to the surface.  The dashed lines indicate the projection of $\tau_{k,*}$'s vertices from the projection point.  When decomposing the volume $\Omega$ the area of the triangle $\tau_{k,*}$, the areas between the arcs on the bounding surface and the edges of $\tau_{k,*}$, and the area of the curved triangle between these arcs make up the boundary of the sliver of volume, $s_{k}$, associated with $t_{k}$.  The projection ensures that between adjacent slivers there are no gaps or overlaps, illustrated by the slivers of volume associated with three adjacent tetrahedra.}
\label{fig:Sliver_Side_Projection}
\end{figure}

When assigning a sliver of volume to a particular tetrahedron, $t_{k}$, care must be taken so that there are no gaps or overlaps between adjacent slivers.  Let $\tau_{k,i}$, $i=1,2,3,4$, be the triangular faces of $t_{k}$.  At least one of these faces has all three vertices on the bounding surface.  In most cases, this will be only one face of $t_{k}$ (particularly when the volume is well resolved by small enough tetrahedra), call it $\tau_{k,*}$.  In the case of the ball, as in \cite{JAR2020}, it turns out that if the three edges of $\tau_{k,*}$ are projected radially from the center of the ball to the surface, gaps and overlaps will be prevented.  Unfortunately this projection does not apply well to every bounded volume.  For each edge of $\tau_{k,*}$ an area between the edge of the triangle and a curve on the bounding surface constructed by projecting the edge to the surface is needed to form a side of the sliver of volume.  The boundary of the sliver volume is formed by all three of these sides, the curved triangle on the bounding surface between the three sides, and the triangle $\tau_{k,*}$.  Figure \ref{fig:Sliver_Side_Projection} illustrates one of these volumes.

To ensure that there are no gaps or overlaps, the procedure presented in \cite{JARBFMLW2016} that is used to partition surfaces with a set of curved triangles is employed here.  This procedure requires first locating a projection point ensuring that for the tetrahedra in the set $\mathcal{K}_{S}$ the triangular faces with three vertices on the bounding surface are not projected onto the same surface area.  Further, the entire bounding surface is included in the union of the projections.

To summarize the definition of the projection point from \cite{JARBFMLW2016}, suppose that the triangle with vertices $\mathbf{a}_{k,*}$, $\mathbf{b}_{k,*}$ and $\mathbf{c}_{k,*}$ is a face of $t_{k}$, and that $\mathbf{a}_{k,*}$, $\mathbf{b}_{k,*}$ and $\mathbf{e}_{k,*}$ are the vertices of a triangular face of an adjacent tetrahedron.  Also suppose that all of these vertices are exactly on the bounding surfaces.  The two triangles share the edge with vertices $\mathbf{a}_{k,*}$ and $\mathbf{b}_{k,*}$.  A unit length normal vector to the triangle with vertices $\mathbf{a}_{k,*}$, $\mathbf{b}_{k,*}$ and $\mathbf{c}_{k,*}$ is defined by
\begin{align}
\mathbf{n}_{\mathbf{a}_{k,*}\mathbf{b}_{k,*}\mathbf{c}_{k,*}}\defeq \frac{(\mathbf{b}_{k,*}-\mathbf{a}_{k,*})\times(\mathbf{c}_{k,*}-\mathbf{a}_{k,*})}{\left\lVert (\mathbf{b}_{k,*}-\mathbf{a}_{k,*})\times(\mathbf{c}_{k,*}-\mathbf{a}_{k,*})\right\rVert_{2}},\nonumber
\end{align}
with $\times$ the vector cross product, and a normal vector for the other triangle can be defined similarly. Here, and in the remainder of this work, the order of the vertices in this definition matters and should be taken as the order shown.

From the normal vectors define
\begin{align}
    \mathbf{u}_{\mathbf{a}_{k,*}\mathbf{b}_{k,*}}=\frac{1}{2}\left(\mathbf{n}_{\mathbf{a}_{k,*}\mathbf{b}_{k,*}\mathbf{c}_{k,*}}+\sign\left(\mathbf{n}_{\mathbf{a}_{k,*}\mathbf{b}_{k,*}\mathbf{c}_{k,*}}^{T}\mathbf{n}_{\mathbf{a}_{k,*}\mathbf{b}_{k,*}\mathbf{e}_{k,*}}\right)\mathbf{n}_{\mathbf{a}_{k,*}\mathbf{b}_{k,*}\mathbf{e}_{k,*}}\right),\nonumber
\end{align}
which is simply an average of the normal vectors to each of these two triangles pointing in roughly the same direction. This vector lies in the plane that separates the slivers of volume assigned to these two triangles.  The plane is also defined to contain the projection point $\mathbf{p}_{k}$ and to be normal to both $\mathbf{u}_{\mathbf{a}_{k,*}\mathbf{b}_{k,*}}$ and the vector $\mathbf{a}_{k,*}-\mathbf{b}_{k,*}$.  The normal vector to this plane that will separate the slivers can then be defined
\begin{align}
\mathbf{v}_{\mathbf{p}_{k}\mathbf{a}_{k,*}\mathbf{b}_{k,*}}=\mathbf{u}_{\mathbf{a}_{k,*}\mathbf{b}_{k,*}}\times(\mathbf{b}_{k,*}-\mathbf{a}_{k,*})\nonumber\
\end{align}
Likewise, the vectors $\mathbf{u}_{\mathbf{b}_{k,*}\mathbf{c}_{k,*}}$, $\mathbf{v}_{\mathbf{p}_{k}\mathbf{b}_{k,*}\mathbf{c}_{k,*}}$, $\mathbf{u}_{\mathbf{c}_{k,*}\mathbf{a}_{k,*}}$ and $\mathbf{v}_{\mathbf{p}_{k}\mathbf{c}_{k,*}\mathbf{a}_{k,*}}$ can be defined from the normal vectors to two other triangular faces of tetrahedra adjacent to $t_{k}$ containing appropriate edges and with all vertices on the bounding surface.

The projection point for each triangle $\tau_{k,*}$ is the intersection of the three planes with normal vectors  $\mathbf{v}_{\mathbf{p}_{k}\mathbf{a}_{k,*}\mathbf{b}_{k,*}}$, $\mathbf{v}_{\mathbf{p}_{k}\mathbf{b}_{k,*}\mathbf{c}_{k,*}}$ and $\mathbf{v}_{\mathbf{p}_{k}\mathbf{c}_{k,*}\mathbf{a}_{k,*}}$ that contain the appropriate vertices/edges.  An expression for this projection point is given by, for instance,
\begin{align}
    \mathbf{p}_{k}=\mathbf{a}_{k,*}+\frac{\mathbf{v}_{\mathbf{p}_{k}\mathbf{b}_{k,*}\mathbf{c}_{k,*}}\cdot(\mathbf{b}_{k,*}-\mathbf{a}_{k,*})}{\mathbf{v}_{\mathbf{p}_{k}\mathbf{b}_{k,*}\mathbf{c}_{k,*}}\cdot\mathbf{w}_{\mathbf{p}_{k}\mathbf{a}_{k,*}}}\mathbf{w}_{\mathbf{p}_{k}\mathbf{a}_{k,*}},\nonumber
\end{align}
where the vector
\begin{align}
\mathbf{w}_{\mathbf{p}_{k}\mathbf{a}_{k,*}}=\mathbf{v}_{\mathbf{p}_{k}\mathbf{a}_{k,*}\mathbf{b}_{k,*}}\times\mathbf{v}_{\mathbf{p}_{k}\mathbf{c}_{k,*}\mathbf{a}_{k,*}},\nonumber
\end{align}
points in the direction of $\mathbf{a}_{k,*}-\mathbf{p}_{k,*}$.

\subsubsection{Evaluating Integrals Over Slivers of Volume Near the Surface When the Implicit Parameterization of the Bounding Surface is Known}  \label{sec:Known_Param}

Assigning the slivers of volume via the projections from $\mathbf{p}_{k}$ provides for a transformation of the coordinates of the sliver which allows the integral over its volume to be written as an iterated integral over a triangular area and a parameter, $\sigma$, which relates to the projection from $\mathbf{p}_{k}$.  Consider any point $\mathbf{x}$ in a neighborhood of the sliver volume.  The line passing through $\mathbf{p}_{k}$ and pointing in the direction of the vector $\mathbf{x}-\mathbf{p}_{k}$ intersects the plane containing $\tau_{k,*}$ at a point $\mathbf{y}_{k,*}$.  The point $\mathbf{x}$ can therefore be represented as
\begin{align}
    \mathbf{x}(\lambda,\mu,\sigma)=\mathbf{y}_{k,*}(\lambda,\mu)+\sigma \mathbf{v}_{k,*}(\lambda,\mu),\nonumber
\end{align}
where
\begin{align}
\mathbf{v}_{k,*}(\lambda,\mu)=\frac{\mathbf{y}_{k,*}(\lambda,\mu)-\mathbf{p}_{k}}{\left\lVert\mathbf{y}_{k,*}(\lambda,\mu)-\mathbf{p}_{k}\right\rVert_{2}}\nonumber \end{align}
and $\lambda$ and $\mu$ define a two dimensional coordinate system for the plane containing $\tau_{k,*}$.  For each value of $\lambda$ and $\mu$ the parameter $\sigma$ ranges between $0$ and $\sigma_{\mbox{max}}(\lambda,\mu)$, where $\sigma_{\mbox{max}}(\lambda,\mu)$ is the smallest value of $\sigma$ such that
$h(\mathbf{x}(\lambda,\mu,\sigma_{\mbox{max}}(\lambda,\mu)))=0$; that is, $\mathbf{x}(\lambda,\mu,\sigma_{\mbox{max}}(\lambda,\mu))$ is on the bounding surface.  The value of $\nu_{k}$ in \eqref{eq:volume_integral_decomposed} is defined as
\begin{align}
\nu_{k}=\sign((\mathbf{m}_{k}-\mathbf{m}_{k,*})\cdot(\mathbf{p}_{k}-\mathbf{m}_{k,*}))\nonumber
\end{align}
to ensure that only integrals over portions of $s_{k}$ outside of $t_{k}$ are added to the total volume and portions of $s_{k}$ inside of $t_{k}$ are subtracted given this choice of coordinate system.  In this definition, $\mathbf{m}_{k}$ and $\mathbf{m}_{k,*}$ are the midpoints (average of the vertices) of the tetrahedron $t_{k}$ and the triangle $\tau_{k,*}$, respectively.

Determining the value of $\sigma_{\mbox{max}}(\lambda,\mu)$ can be simple if the function $h$ that defines the bounding surface is known and has roots that can be written in closed form.  If the $h$ is known but the roots cannot be written in closed form, then a rootfinding method such as Newton's method or Broyden's method \cite{CGB1965} must be used to find the root, $\sigma_{\mbox{max}}(\lambda,\mu)$, of $h(\mathbf{x}(\lambda,\mu,\sigma_{\mbox{max}}(\lambda,\mu)))=0$ for each value of $\lambda$ and $\mu$ required to approximate the integral over the sliver volume.  In either case, when $h$ is known the parameterization of $\tau_{k,*}$ via, for instance,
\begin{align}
\mathbf{y}_{k,*}(\lambda,\mu)=(1-\lambda) \mathbf{a}_{k,*}+\lambda((1-\mu)\mathbf{b}_{k,*}+\mu\mathbf{c}_{k,*})\nonumber
\end{align}
where $0\leq\lambda\leq 1$ and $0\leq\mu\leq1$, allows the use of any number of quadrature techniques over intervals to evaluate the iterated integrals in
\begin{align}
    \iiint\limits_{s_{k}}\phi\left(\left\lVert \mathbf{x}-\mathbf{x}_{k,j}\right\rVert\right)dV=\int\limits_{0}^{1}\int\limits_{0}^{1}\int\limits_{0}^{\sigma_{\mbox{max}}(\lambda,\mu)}\phi\left(\left\lVert \mathbf{x}(\lambda,\mu,\sigma)-\mathbf{x}_{k,j}\right\rVert\right) \left\lvert J(\lambda,\mu,\sigma)\right\rvert d\sigma d\lambda d\mu. \nonumber
\end{align}
In the right hand side,
\begin{align}
J(\lambda,\mu,\sigma)=\left(1+\frac{\sigma}{\lVert \mathbf{y}_{k,*}(\lambda,\mu)-\mathbf{p}_{k}\rVert_{2}}\right)^{2}\lambda\left[ \mathbf{v}_{k,*}(\lambda,\mu)\cdot((\mathbf{b}_{k,*}-\mathbf{a}_{k,*})\times(\mathbf{c}_{k,*}-\mathbf{b}_{k,*}))\right]\nonumber
\end{align}
is the Jacobian determinant of $\mathbf{x}$ with respect to $\sigma$, $\lambda$ and $\mu$.   On the other hand, if $h$ is not known, a procedure for computing $\sigma_{\mbox{max}}(\lambda,\mu)$ and quadrature weights for evaluating the integrals in $\lambda$ and $\mu$ is described in the following section.  This procedure takes advantage of the knowledge that there are points in $\mathcal{S}_{N}$ near $\tau_{k,*}$, including its vertices, that are also on the bounding surface.  Once $\sigma_{\mbox{max}}(\lambda,\mu)$ is known at appropriate values of $\lambda$ and $\mu$, the integral in $\sigma$ can be easily treated with any number of quadrature methods over an interval.  In the current implementation, a 21 node Legendre-Gauss-Lobatto quadrature rule is used for evaluating the integrals in each of $\lambda$, $\mu$ and $\sigma$.

\subsubsection{Evaluating Integrals Over Slivers of Volume Near the Surface When the Implicit Parameterization of the Bounding Surface is Unknown} \label{sec:Unknown_Param}

When $h$ is unknown or unavailable, for instance, when the node set $\mathcal{S}_{N}$ is determined by simulation or physical measurement, the value of $\sigma_{\mbox{max}}(\lambda,\mu)$ cannot be computed for most choices of $\lambda$ and $\mu$.  Therefore, a procedure is needed that allows the approximation of the integral over $s_{k}$  using only points from $\mathcal{S}_{N}$ that are located on the bounding surface where the upper bound on $\sigma$ is directly computable.

Define a two dimensional coordinate system for the plane containing $\tau_{k,*}$ via
\begin{align}
\mathbf{y}_{k,*}(\lambda,\mu)=R_{k}^{T}\left[\begin{array}{c}\lambda \\ \mu \\ g_{k}\end{array}\right]+\mathbf{p}_{k},\nonumber
\end{align}
where
\begin{align}
g_{k}=\lvert(\mathbf{m}_{k,*}-\mathbf{p}_{k})\cdot \mathbf{n}_{\mathbf{a}_{k,*}\mathbf{b}_{k,*}\mathbf{c}_{k,*}}\rvert,\nonumber
\end{align}
\begin{align}
\mathbf{m}_{k,*}=\frac{1}{3}(\mathbf{a}_{k,*}+\mathbf{b}_{k,*}+\mathbf{c}_{k,*})\nonumber
\end{align}
and
\begin{align}
R_{k}=& \sign(n_{k,x})\sign(n_{k,z}) \nonumber\\
&\left\lbrack\begin{array}{ccc}
\frac{n_{k,z}n_{k,x}}{\sqrt{n_{k,x}^{2}+n_{k,y}^{2}}} & \frac{n_{k,z}n_{k,y}}{\sqrt{n_{k,x}^{2}+n_{k,y}^{2}}} & -\frac{(n_{k,x}^{2}+n_{k,y}^{2})}{\sqrt{n_{k,x}^{2}+n_{k,y}^{2}}} \\
-\frac{n_{k,y}\sign(n_{k,z})}{\sqrt{n_{k,x}^{2}+n_{k,y}^{2}}} & \frac{n_{k,x}\sign(n_{k,z})}{\sqrt{n_{k,x}^{2}+n_{k,y}^{2}}} & 0 \\
\lvert n_{k,x}\rvert & n_{k,y}\sign(n_{k,x}) & n_{k,z}\sign(n_{k,x})
\end{array}\right\rbrack,\nonumber
\end{align}
with $n_{k,x}$, $n_{k,y}$ and $n_{k,z}$ are the three components of $\mathbf{n}_{\mathbf{a}_{k,*}\mathbf{b}_{k,*}\mathbf{c}_{k,*}}$.

Let $\mathcal{N}_{k,*}=\left\{\mathbf{x}_{k,*,j}\right\}_{j=1}^{\eta}$ be the $\eta$ points in $\mathcal{S}_{N}$ nearest to the midpoint, $\mathbf{m}_{k,*}$, of $\tau_{k,*}$ that also satisfy $h(\mathbf{x}_{k,*,j})=0$.  These are points exactly on the bounding surface. There is a point at the intersection of the plane containing $\tau_{k,*}$ and the line segment connecting $\mathbf{x}_{k,*,j}$ and $\mathbf{p}_{k}$ defined by
\begin{align}
    \mathbf{y}_{k,*,j}=\mathbf{x}_{k,*,j}-\frac{(\mathbf{x}_{k,*,j}-\mathbf{m}_{k})\cdot\mathbf{n}_{\mathbf{a}_{k,*}\mathbf{b}_{k,*}\mathbf{c}_{k,*}}}{(\mathbf{x}_{k,*,j}-\mathbf{p}_{k})\cdot\mathbf{n}_{\mathbf{a}_{k,*}\mathbf{b}_{k,*}\mathbf{c}_{k,*}}}(\mathbf{x}_{k,*,j}-\mathbf{p}_{k}).\nonumber
\end{align}
Denote the values of $\lambda$ and $\mu$ at these intersection points as $\lambda_{k,j}$ and $\mu_{k,j}$, that is $\mathbf{y}_{k,*}(\lambda_{k,j},\mu_{k,j})=\mathbf{y}_{k,*,j}$ (and $\mathbf{v}_{k,*}(\lambda_{k,j},\mu_{k,j})=\mathbf{v}_{k,*,j}$), so that
\begin{align}
\sigma_{\mbox{max}}(\lambda_{k,j},\mu_{k,j})=\frac{(\mathbf{x}_{k,*,j}-\mathbf{y}_{k,*,j})\cdot\mathbf{n}_{\mathbf{a}_{k,*}\mathbf{b}_{k,*}\mathbf{c}_{k,*}}}{\mathbf{v}_{k,*,j}\cdot\mathbf{n}_{\mathbf{a}_{k,*}\mathbf{b}_{k,*}\mathbf{c}_{k,*}}}.\nonumber
\end{align}
Now the integral over the sliver volume becomes
\begin{align}
    \iiint\limits_{s_{k}}\phi\left(\left\lVert \mathbf{x}-\mathbf{x}_{k,j}\right\rVert\right)dV=\iint\limits_{\tau_{k,*}}\int\limits_{0}^{\sigma_{\mbox{max}}(\lambda,\mu)}\phi\left(\left\lVert \mathbf{x}(\lambda,\mu,\sigma)-\mathbf{x}_{k,j}\right\rVert\right) \left\lvert J(\lambda,\mu,\sigma)\right\rvert d\sigma d\lambda d\mu, \nonumber
\end{align}
with the Jacobian determinant now
\begin{align}
J(\lambda,\mu,\sigma)=\left(1+\frac{\sigma}{\lVert \mathbf{y}_{k,*}(\lambda,\mu)-\mathbf{p}_{k}\rVert_{2}}\right)^{2}\left[ \mathbf{v}_{k,*}(\lambda,\mu)\cdot\mathbf{n}_{\mathbf{a}_{k,*}\mathbf{b}_{k,*}\mathbf{c}_{k,*}}\right].\nonumber
\end{align}
In this expression the quadrature weights for approximating the integral over $\tau_{k,*}$ can be found by now constructing a two dimensional RBF interpolant over the set of points $\left\{(\lambda_{k,j},\mu_{k,j})\right\}_{j=1}^{\eta}$ and integrating the interpolant, similar to what is presented in sections \ref{sec:ConstructInterpolant} and \ref{sec:IntegrateInterpolant}.  In this case, the system of linear equations will require the RBFs and any bivariate polynomial terms included in the basis set to be integrated over the projection of $\tau_{k,*}$ into the $(\lambda,\mu)$ coordinate system.  Closed form expressions for the integrals of the polynomial terms are easy to compute, and for many common RBFs closed form expressions also exist and can be computed as in \cite{JARBF2016,JARBFMLW2016,JARBF2017}.  The integral over $\sigma$ can be approximated using any number of quadrature methods for an interval.  In the current implementation, a 21 node Legendre-Gauss-Lobatto quadrature rule is used for evaluating the integral in $\sigma$.

\subsection{Step 4: Combine weights from the subdomains}
Once the weights for each subdomain are computed, summing over all $k\in\{1,2,\ldots,K\}$ leads to the approximation of the volume integral over $\Omega$
\begin{align}
\iiint\limits_{\Omega}f(\mathbf{x})dV\approx\sum\limits_{k=1}^{K}\sum_{j=1}^{n}w_{k,j}f(\mathbf{x}_{k,j}). \nonumber
\end{align}

Let $\mathcal{K}_{i}$, $i=1,2,\ldots,N$, be the set of all pairs $(k,j)$ such that $\mathbf{x}_{k,j}\mapsto\mathbf{x}_{i}$. Then the volume integral over $\Omega$ can be rewritten as
\begin{align}
\iiint\limits_{\Omega}f(\mathbf{x})dV\approx\sum\limits_{i=1}^{N}W_{i}f(\mathbf{x}_{i}).\nonumber
\end{align}

\section{Test Examples} \label{sec:Test_Examples}

\begin{figure}[h]
\begin{center}
\includegraphics[width=\linewidth]{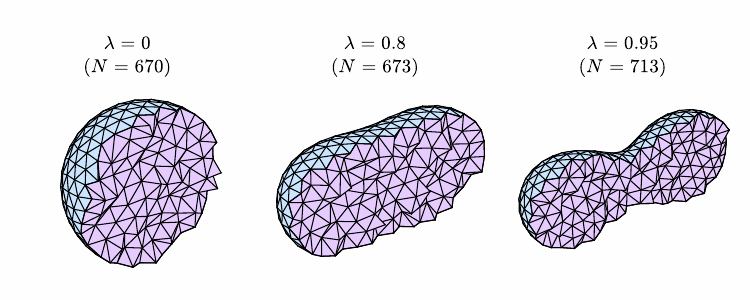}
\end{center}
\caption{Examples of the quasi-unifomly spaced nodes in the volumes with bounding surfaces defined implicitly by \eqref{eq:surface_parameterization}.}
\label{fig:Bounded_Volume_Node_Sets}
\end{figure}

To demonstrate the performance of the method describe in section \ref{sec:Algorithm_Description}, the algorithm was applied to three different test integrands featuring varying degrees of smoothness.  Similar to the tests described in \cite{JARBFMLW2016} and \cite{JARBF2017}, weights are computed on quasi-uniformly spaced nodes bounded by a family of surfaces generated by rotating the Cassini ovals about the $x$-axis.  That is, the bounding surfaces for these volumes are defined by the level surfaces
\begin{align}
    h(\mathbf{x})=h(x,y,z)=(x^{2}+y^{2}+z^{2})^{2}-2\alpha^{2}(x^{2}-y^{2}-z^{2})+\alpha^{4}-\beta^{4}=0,\label{eq:surface_parameterization}
\end{align}
which depend on the two parameters $\alpha$ and $\beta$. The demonstrations here will consider $\alpha=\lambda \beta$ for $0<\lambda<1$ (specifically, $\lambda=0,0.8,0.95$).  The parameter $\beta$ in this work is chosen so that the volume is equal to 1, with $\beta$ chosen numerically.  Quadrature nodes for these volumes were generated using a modification of the algorithm presented in \cite{PPGS04}. Examples of the node sets are displayed in figure \ref{fig:Bounded_Volume_Node_Sets}.

\subsection{Performance on Test Integrands}

The algorithm was applied to three test integrands, first by utilizing the knowledge of the surface parameterization as in section \ref{sec:Known_Param} and second disregarding the known surface parameterization as in section \ref{sec:Unknown_Param}.  Just as in \cite{JAR2020}, when generating quadrature weights the radial basis function when constructing the interpolant on the \textit{volume} was $\phi(r)=r^{3}$ and the number of nearest neighbors, $n=\frac{(m+1)(m+2)(m+3)}{3}$, was based on the trivariate polynomial order $m$. The choice for $n$ was guided by some computational experiments in, for instance, \cite{JARBF2017,VBNFBFGAB2017,VBNFBF2019} that indicated that in the presence of boundaries the number of nearest neighbors must be large enough to overcome effects like Runge phenomenon.  The examples given in \cite{JARBF2017} indicated that the boundary errors were most prominent when nodes were (exactly) uniformly spaced, and the experiments in \cite{JAR2020} suggested that this is a safe choice. When considering the case of an unknown bounding surface parameterization, the RBF chosen for integrals in the \textit{plane} was $r^{7}$ and the number of nearest neighbors was $\eta=\left\lceil1.05\frac{(\gamma+1)(\gamma+2)}{2}\right\rceil$. Here $\gamma$ is the order of the bivariate polynomial included in the basis and in all cases $\gamma=2m$. Given that the triangle to be integrated over is on the interior of the set of nodes used for interpolation, boundary effects pose less of a problem and the choice of $\eta$ appears safe in this context.

\begin{figure}[h]
\begin{center}
\includegraphics[width=\linewidth]{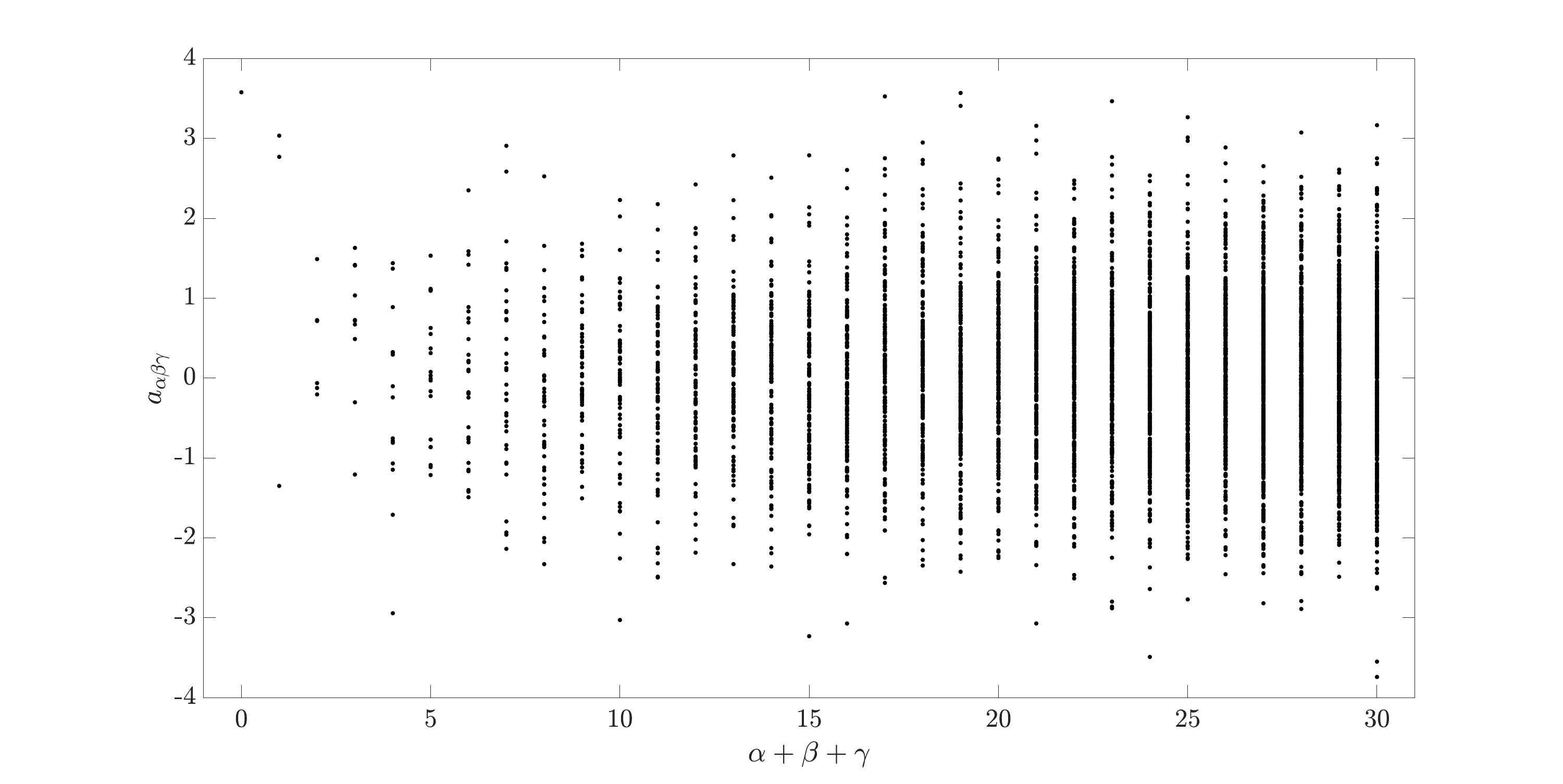}
\end{center}
\caption{Values of the coefficients $a_{\alpha\beta\gamma}$ in $f_{1}(x,y,z)$ shown relative to total order $\alpha+\beta+\gamma$ of the term that they are multiplying.}
\label{fig:f1_Coefficients}
\end{figure}

\begin{figure}[h]
\begin{center}
\includegraphics[width=\linewidth]{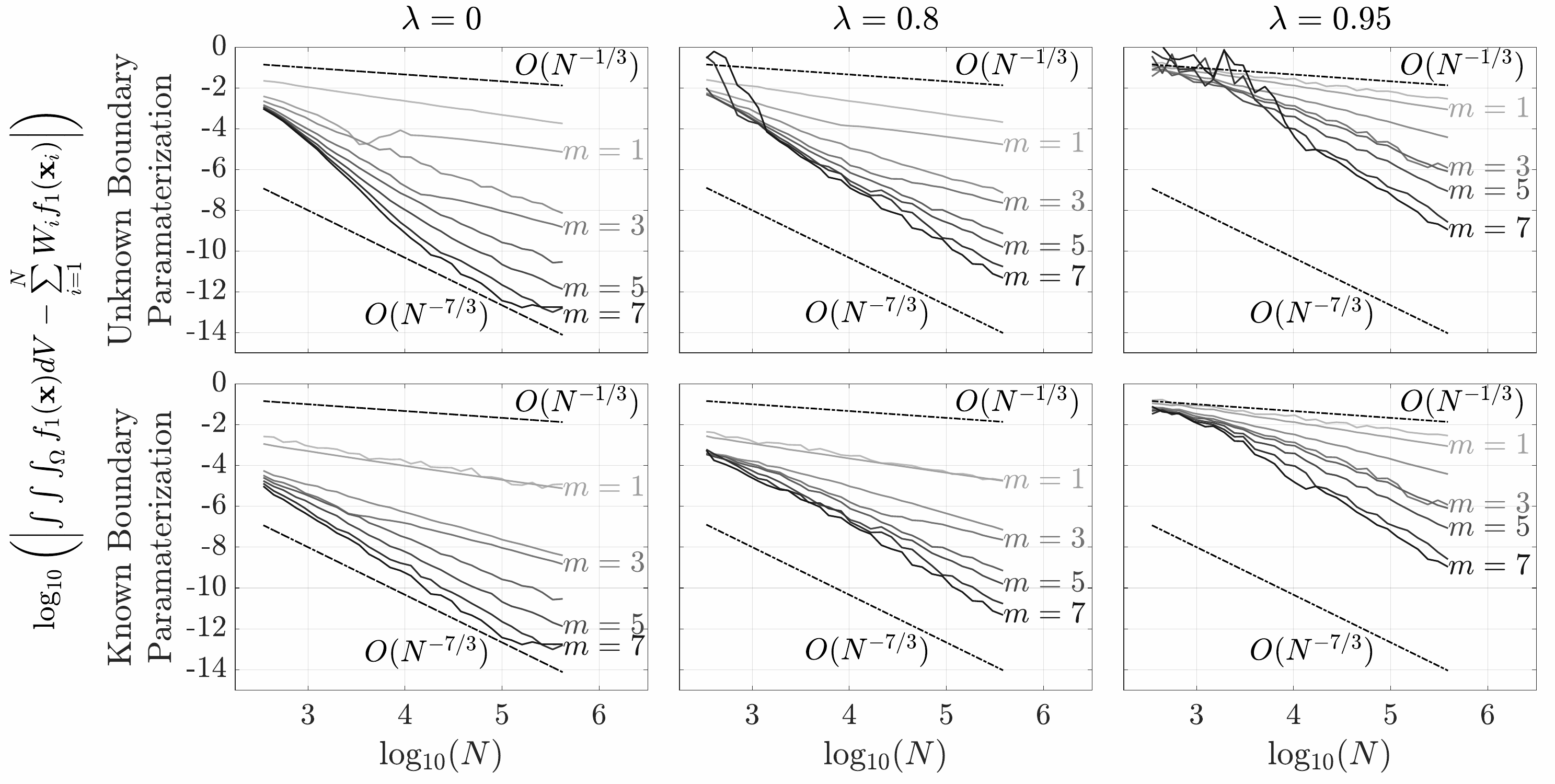}
\end{center}
\caption{Log base 10 of the absolute error when approximating the volume integral of $f_1$ over the test surfaces.  The errors shown here are the largest after rotating the integrand 1000 times.  The shade of the curves darken as the total order, $m$, of the trivariate polynomials that are included in the basis set increases (with only odd orders labeled in the figure).}
\label{fig:Bounded_Volume_Quadrature_Error_f1}
\end{figure}

The first of the test integrands is a degree 30 trivariate polynomial.  That is, let
\begin{align}
f_{1}(x,y,z)=\sum_{\alpha+\beta+\gamma=0}^{30}a_{\alpha\beta\gamma}x^{\alpha}y^{\beta}z^{\gamma},\nonumber
\end{align}
with $\alpha\geq0$, $\beta\geq0$ and $\gamma\geq0$.  Figure \ref{fig:f1_Coefficients} displays the values of the coefficients $a_{\alpha\beta\gamma}$ (exact values are available at \cite{VolumeQuadCodeMatlab}) relative to total order $\alpha+\beta+\gamma$ of the term that they are multiplying.  From this we can see that at each order there are coefficients of $O(1)$ so that for $\lVert\mathbf{x}\rVert_{2}$ small enough, this sum has terms that decay like $O(\lVert\mathbf{x}\rVert_{2}^{\alpha+\beta+\gamma})$ until $\alpha+\beta+\gamma=30$, beyond which the terms are zero.
Figure \ref{fig:Bounded_Volume_Quadrature_Error_f1} displays convergence of the approximate integral to the exact value at an order better than O$\left(N^{-\frac{m}{3}}\right)$, where $m$ corresponds to the order of the trivariate polynomial terms used in the approximation.  Although only $O(N^{-\frac{1}{3}})$ and $O(N^{-\frac{7}{3}})$ reference curves are shown, if $\Delta$ refers to a typical node separation distance, this corresponds to a convergence order of better than O$\left(\Delta^{m}\right)$ at each order. The theory in \cite{VB2019} explains that if the multivariate polynomial basis up to degree $m$ is included in the process of RBF interpolation, then all of the terms in a series expansion up to degree $m$ will be handled exactly for the function being interpolated.  The remaining terms in the series are then approximated by the RBF basis that was included, leading to convergence at the order of at least O$(\Delta^{m})$, regardless of whether the boundary surface parameterization is known.  In the cases where the boundary surface parameterization is unknown the errors are larger for $N$ less than between $10^4$ and $10^5$, particularly in the case of $\lambda=0.95$.  Consistent with the results in \cite{JARBFMLW2016}, which inspired the method for integrating over sliver volumes near the surface when the parameterization is unknown, these larger errors are the result of the projection from $\mathbf{p}_{k}$ to the plane containing $\tau_{k,*}$ sending points that are relatively close (in Euclidean distance) on the boundary surface to points potentially far away in the plane when the direction of the surface normal changes rapidly locally.

%

\begin{figure}[h]
\begin{center}
\includegraphics[width=\linewidth]{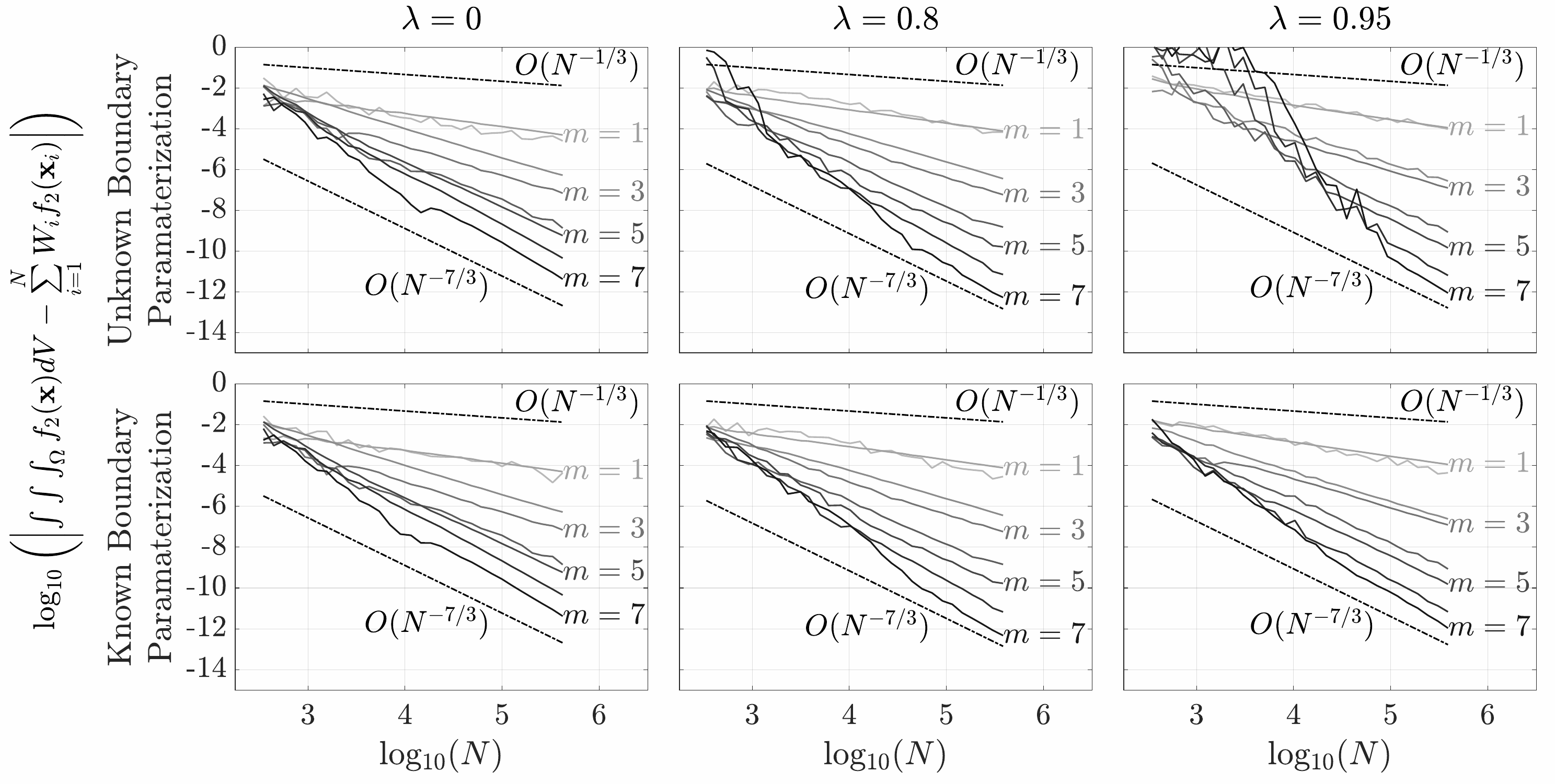}
\end{center}
\caption{Log base 10 of the absolute error when approximating the volume integral of $f_2$ over the test surfaces.  The errors shown here are the largest after rotating the integrand 1000 times.  The shade of the curves darken as the total order, $m$, of the trivariate polynomials that are included in the basis set increases (with only odd orders labeled in the figure).}
\label{fig:Bounded_Volume_Quadrature_Error_f2}
\end{figure}

The second test integrand is the Gaussian
\begin{align}
    f_{2}(x,y,z)=\exp\left(-10\left((x-x_{s})^{2}+(y-y_{s})^{2}+(z-z_{s})^2\right)\right)\nonumber
\end{align}
where
\begin{align}
(x_{s},y_{s},z_{s})=(0.047056440432708,0.071766893999009,0.118950756342700)\nonumber
\end{align}
is a randomly chosen shift of the center of the Gaussian from the origin. In order to have an accurate value to compare to, the volume integral was first approximated by evaluating
\begin{align}
\int\limits_{-\beta\sqrt{1+\lambda^2}}^{\beta\sqrt{1+\lambda^2}}\int\limits_{-\sqrt{\sqrt{\beta^4 + 4 \beta^2 \lambda^2 x^2} - x^2 - \beta^2 \lambda^2}}^{\sqrt{\sqrt{\beta^4 + 4 \beta^2 \lambda^2 x^2} - x^2 - \beta^2 \lambda^2}}\int\limits_{-\sqrt{\sqrt{\beta^4 + 4 \beta^2 \lambda^2 x^2} - x^2 -y^2- \beta^2 \lambda^2}}^{\sqrt{\sqrt{\beta^4 + 4 \beta^2 \lambda^2 x^2} - x^2 - y^2 - \beta^2 \lambda^2}}f_{2}(x,y,z)dzdydx,\nonumber
\end{align}
using Matlab's integral3 command with the absolute and relative tolerances both set to ten times machine precision.  Figure \ref{fig:Bounded_Volume_Quadrature_Error_f2} illustrates the error in the integral of $f_{2}$ over each volume when compared to the result from Matlab after rotating the integrand randomly 1000 times.  It is clear again that the order of the error is most dependent on the degree of the polynomials used in the interpolation.  This can be explained by considering that the series
\begin{align}
    f_{2}(x,y,z)=\sum\limits_{\zeta_{x}=0}^{\infty}\sum\limits_{\zeta_{y}=0}^{\infty}\sum\limits_{\zeta_{z}=0}^{\infty}\frac{(-10)^{\zeta_{x}+\zeta_{y}+\zeta_{z}}}{(\zeta_{x}!)(\zeta_{y}!)(\zeta_{z}!)}(x-x_{s})^{2\zeta_{x}}(y-y_{s})^{2\zeta_{y}}(z-z_{s})^{2\zeta_{z}}\nonumber
\end{align}
has terms that decay rapidly as the total order $2(\zeta_{x}+\zeta_{y}+\zeta_{z})$ increases so that after all terms of degree $m$ are accounted for exactly, only terms of size $O(\lVert\mathbf{x}-\mathbf{x}_{s}\rVert_{2}^{m+1})$ or smaller remain.

\begin{figure}[h]
\begin{center}
\includegraphics[width=\linewidth]{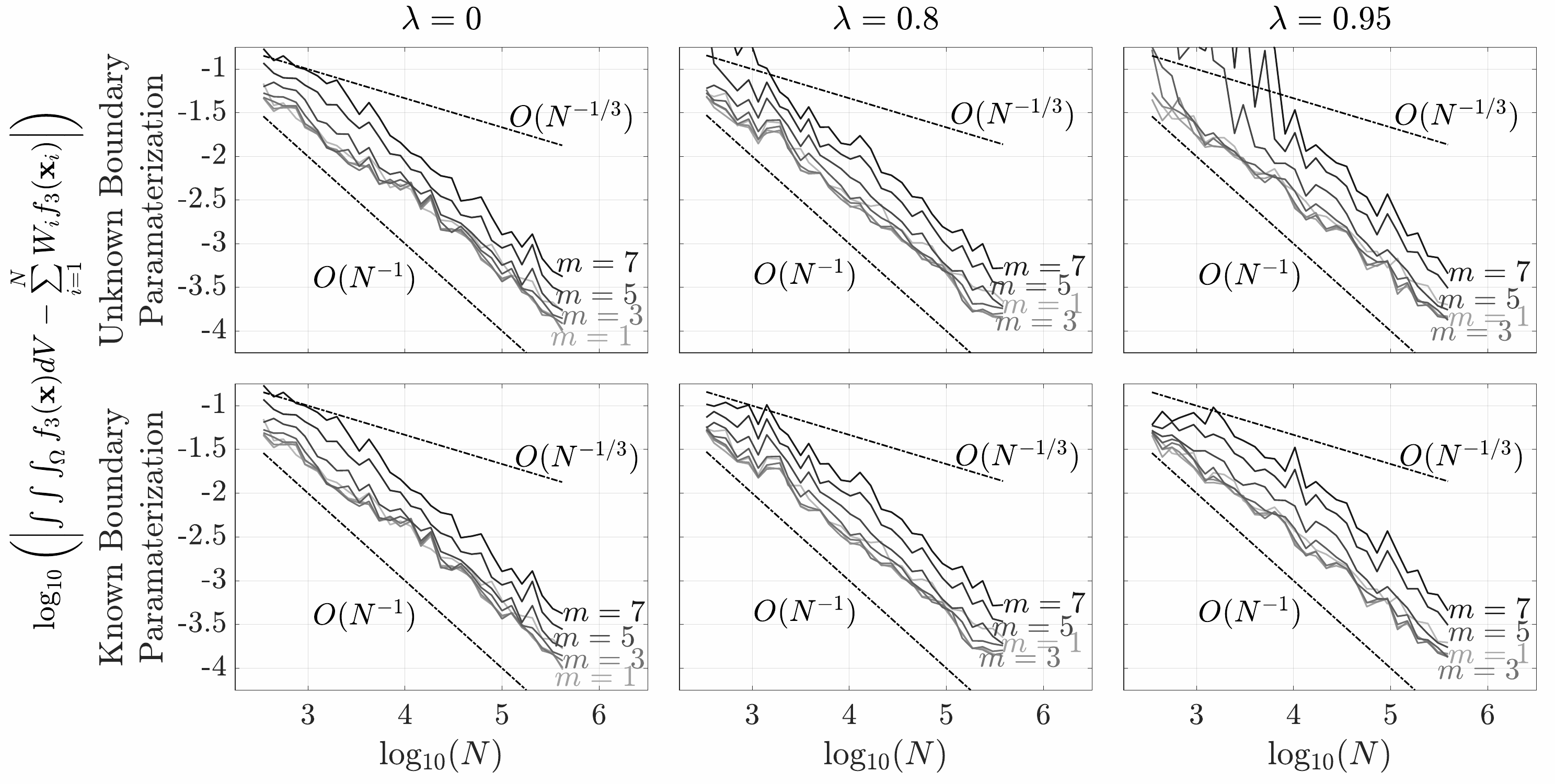}
\end{center}
\caption{Log base 10 of the absolute error when approximating the volume integral of $f_3$ over the test surfaces.  The errors shown here are the largest after rotating the integrand 1000 times.  The shade of the curves darken as the total order, $m$, of the trivariate polynomials that are included in the basis set increases (with only odd orders labeled in the figure).}
\label{fig:Bounded_Volume_Quadrature_Error_f3}
\end{figure}

Finally, the algorithm was applied to a test integrand featuring a steep and localized gradient.  The test integrand was
\begin{align}
    f_{3}(x,y,z)=\tan^{-1}\left(500 z\right),\nonumber
\end{align}
which has a steep gradient near the plane $z=0$.  The power series for this function has the form
\begin{align}
f_{3}(x,y,z) = \sum\limits_{\zeta=1}^{\infty}(-1)^{\zeta-1}\frac{500^{2\zeta-1}}{2\zeta-1}z^{2\zeta-1},\nonumber
\end{align}
which has terms, for $\zeta$ large enough, that grow for any fixed choice of $z>0$.  Therefore, even when capturing the terms of the power series up to order $m$ exactly (by including trivariate polynomials up to order $m$ in \eqref{eq:rbf_interpolant}), the accuracy of the approximation depends on how well the RBF basis can approximate the remainder of the power series.  Figure \ref{fig:Bounded_Volume_Quadrature_Error_f3} illustrates this with the error decaying at roughly the same rate (slightly less than $O(\Delta^{3})$ since $\Delta\sim N^{-1/3}$) regardless of the included polynomial order.  The numerical tests in \cite{JAR2020} suggest that node sets capturing the rapid change in the integrand can significantly improve the performance of the algorithm.

\subsection{Computational Expense}

\begin{figure}[h]
\begin{center}
\includegraphics[width=0.66\linewidth]{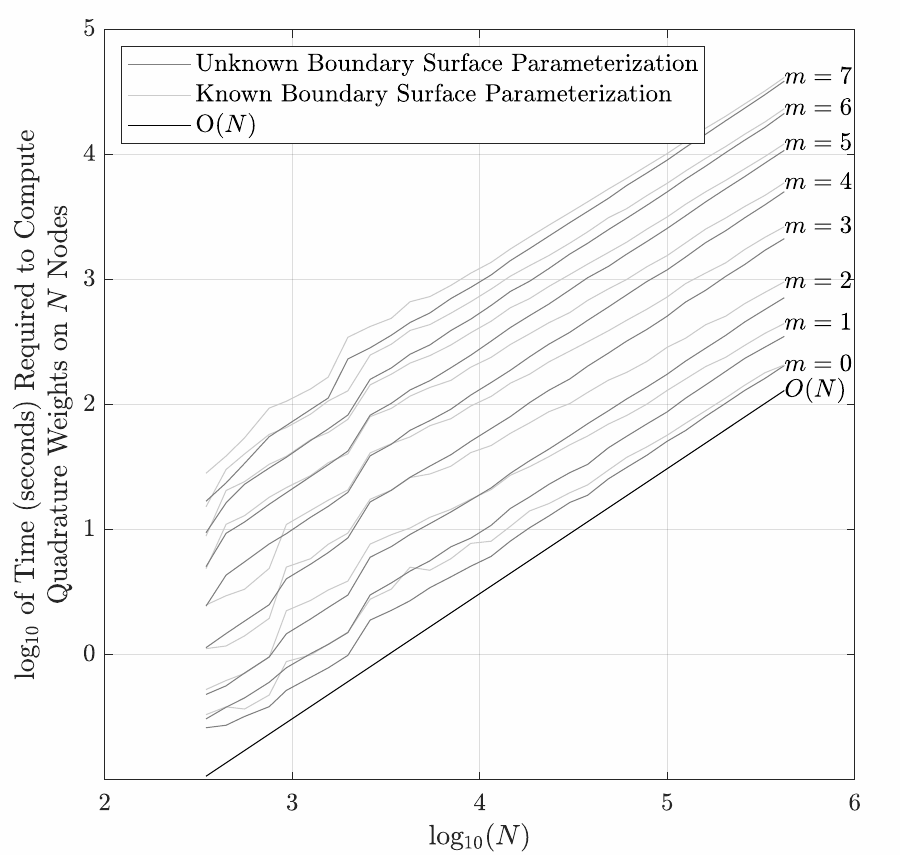}
\end{center}
\caption{Log base 10 of the time it takes to compute quadrature weights on $N$ nodes when including trivariate polynomial up to order $m$.  }
\label{fig:Bounded_Volume_Quadrature_Times}
\end{figure}

Figure \ref{fig:Bounded_Volume_Quadrature_Times} illustrates the time to compute the set of quadrature weights on $N$ nodes for various choice of the polynomial order, $m$.  Consideration of each tetrahedron individually allows the time to compute a set of quadrature weights (and the use of memory) to scale like $O(N)$.  This is in line with what is presented in \cite{JAR2020,JARBF2016,JARBFMLW2016,JARBF2017}.  Since the choice of $m$ (and $\mu$) affects the sizes of the systems of linear equations that need to be solved at each iteration, the figure shows an increase in the computational cost as $m$ increases.  That is, for each tetrahedron a system of linear equations of size $n+M$ must be solved at a cost of $O((n+M)^{3})$, with $M$ and $n$ both dependent on $m$.  The cost of computing quadrature weights for each tetrahedron are also dependent on parameter choices for the quadrature rules used in sections \ref{sec:Known_Param} and \ref{sec:Unknown_Param} when integrating the three-dimensional RBF basis over the sliver volumes; however, the overall cost still scales as $O(N)$.
Although parallelization tests are not presented for this specific algorithm, the results presented in \cite{JARBF2016} should translate to the algorithm proposed here.  Except for the identification of nearest neighbors in order to construct the local weight set for each tetrahedron/sliver of volume and for the combination of weights in step 4 the algorithm is pleasingly parallel.

\section{Conclusions} \label{sec:Conclusions}

This study has augmented the previous RBF-FD based approach for evaluating definite integrals \cite{JARBF2016,JARBFMLW2016,JARBF2017,JAR2020} with an extension to integrals over volumes bounded by smooth surfaces. An important aspect of this extension is that explicit knowledge of an expression for the bounding surface need not be known.  While the tests were performed on quasi-uniformly spaced node sets, spatially varying density of the nodes is permitted and can be leveraged when the integrand or volume of integration require finer resolution in certain places as illustrated in \cite{JAR2020}.  The computational tests illustrate an algorithm that can achieve at least $O(\Delta^{m})$ accuracy, with $\Delta$ the typical node separation distance and $m$ the order of trivariate polynomial basis functions included in the approximation.  On a set of $N$ nodes in the volume, the computational cost is only $O(N)$ and the algorithm is pleasingly parallel.

\bibliographystyle{unsrt}      
\bibliography{references}   


\end{document}